\documentclass[twoside,11pt]{article}
\usepackage[section]{algorithm}
\usepackage{amsthm,amsmath,amssymb,amsmath,amsfonts,graphicx,fancyhdr,algorithmic,longtable,listings,color,xspace}
\usepackage{bibcheck}
\usepackage{epstopdf}
\def\enddemo{\qed \endtrivlist}
\expandafter\let\csname enddemo*\endcsname=\enddemo

\def\qedsymbol{\ifmmode\bgroup\else$\bgroup\aftergroup$\fi
  \vcenter{\hrule\hbox{\vrule
height.6em\kern.6em\vrule}\hrule}\egroup}
\def\qed{\ifmmode\else\unskip\nobreak\fi\quad\qedsymbol}

\newtheorem{thm}{Theorem}[section]

\theoremstyle{definition}

\newtheorem{rem}{Remark}[section]
\theoremstyle{plain}

\newcommand{\bb}{\begin{equation}}
\newcommand{\ee}{\end{equation}}

\newcommand{\R}{{\mathbb R}}

\def\mx{\mbox{\matbi{X}}}
\def\mc{\mbox{\matbi{C}}}
\def\my{\mbox{\matbi{Y}}}

\def\mcs{\mbox{\matbi{c}}}

\def\da{\d{${A}$}}
\def\dh{\d{${H}$}}
\def\dx{\d{${{X}}$}}
\def\di{\d{${{I}}$}}
\def\dy{\d{${{Y}}$}}
\def\ds{\d{${{S}}$}}
\def\dt{\d{${{T}}$}}
\def\dm{\d{${{M}}$}}
\def\dq{\d{${{Q}}$}}
\def\dr{\d{${{R}}$}}
\def\dc{\d{${{C}}$}}

\def\ba{\begin{eqnarray*}}
\def\ea{\end{eqnarray*}}
\def\baa{\begin{eqnarray}}
\def\eaa{\end{eqnarray}}

\def\tfrac{\textstyle\frac}

\font\matbi=cmbxti10

\def\la{\lambda}
\def\cl{\centerline}
\def\od{{\cal{O}}}

\font\rm=cmss11

\font\bsfs=cmssbx10 at 14.4pt

\font\bsfv=cmssbx10 at 17.28pt

\title{\bsfv Iterative methods for the inclusion of the inverse matrix}
\begin{footnotesize}
\author{\frenchspacing
{\bf Marko D. Petkovi\' c}\footnote{Corresponding author. Emails: dexterofnis@gmail.com (M.D. Petkovi\'c),  miodragpetkovic@gmail.com (M.S. Petkovi\'c).}\\
{\small \it University of Ni\v{s}, Faculty of Science and Mathematics}\\
{\small \it Vi\v segradska 33, 18000 Ni\v s, Serbia}
 \\[8pt]
{\bf Miodrag S. Petkovi\' c} \\
{\small \it University of Ni\v{s}, Faculty of Electronic Engineering}\\
{\small \it Aleksandra Medvedeva 14, 18000 Ni\v s, Serbia} \\
}

\end{footnotesize}
\date{}
\begin{document}

\maketitle

\begin{abstract} {\fontfamily{cmss}\selectfont \fontseries{msc}
\selectfont
\fontshape{n}\selectfont \fontsize{9}{9}\selectfont
In this paper we present an efficient iterative method of order six for the inclusion of the inverse of a given regular matrix. To provide the upper error bound of the outer matrix for
the inverse matrix, we  combine point and interval iterations. The new method is relied on a suitable matrix identity and a modification of a hyper-power method. This method is also feasible in the case of a full-rank $m\times n$ matrix,  producing the interval sequence which converges to the Moore-Penrose inverse. It is shown that computational efficiency of the proposed method is equal or higher than the methods of hyper-power's type.}

\frenchspacing \itemsep=-1pt
\begin{description}
\item[] {\it AMS Subject Classification}: 15A09, 65G30, 47J25, 03D15, 65H05.
\item[] {\it Key words}: Inclusion methods; inverse matrix; hyper-power methods; convergence; computational efficiency.
\end{description}
\end{abstract}

\renewcommand{\thefootnote}{}%
\footnote{This work is supported by the Serbian
Ministry of Education and Science under the grants 174033 (first author) and 174022 (second author).}

{
\fontseries{msc}
\fontsize{11}{12}\selectfont

\section{\bsfs Introduction}
\setcounter{equation}{0}

A number of tasks in Numerical analysis, Graph theory, Geometry, Statistics, Computer sciences, Cryptography (encoding and decoding matrices), Partial differential equations, Physics, Engineering disciplines, Medicine (eg., digital tomosynthesis), Management and Optimization (Design Structure Matrix) and so on, is modeled in the matrix form. Solution of these problems is very often reduced to finding an inverse matrix. There is a vast literature in this area so that we will not consider all matrix numerical methods of iterative nature. Instead, in this paper we concentrate only on that small branch of matrix iterative analysis concerned with the efficient determination of inverse matrices
 with upper error bound of the solution using interval arithmetic.  The presented study is a two-way bridge between linear algebra and computing.

The paper is divided into four sections and organized as follows.  In Section 2
we give some preliminary matrix properties and definitions and a short study of hyper-power matrix iterations. The main goal of this paper is to state an efficient iterative method of order six for the inclusion of the inverse of a given regular matrix, which is the subject of Section 3. This method is constructed by modifying a hyper-power method in such a way that the computational cost is decreased. In order to provide information on the upper error bounds of the approximate interval matrix,  interval arithmetic is used.
  Computational aspects of the considered interval methods and one numerical example are considered in Section 4. We show that computational efficiency of the proposed method is equal or higher than the methods of hyper-power's type realized in a Horner scheme fashion.

\section{\bsfs Hyper-power methods}

Applying  numerical methods on digital computers, one of the most important task is to provide an information on the accuracy of obtained results. The interest of bounding roundoff errors in matrix computations has come from the impossibility of exact representation of elements of matrices in some cases since numbers are represented in the computer by string of bits of fixed, {\it finite} length. For more details see  \cite{mur3}, \cite{herz}, \cite{cia}. Such case also appears in finding inverse matrices, the subject of this paper. To provide the upper error bound of the outer matrix for
the inverse matrix, we will combine point and interval iterations.
The essential advantage  of the presented interval methods consists of  capturing all the roundoff errors automatically, making this approach useful, elegant and powerful tool
for finding errors in the sought results.

To avoid any confusion, in this paper interval matrices will be
denoted by bold capital letters and real matrices (often called
point matrices) by calligraphic letters with a dot below  the
letter. We use bold small letters to denote real intervals.

Let $\dc=[\mcs
_{ij}]$ be a nonsingular $n\times n$ matrix, where
$\mcs_{ij}=[\underline{c}_{\,ij},\overline{c}_{ij}],\
\overline{c}_{ij}-\underline{c}_{\,ij}\ge 0,$ are real intervals.
 An interval matrix whose all elements are points (real numbers) is called a {\it point matrix}. Basic definitions, operations and properties of interval matrices can be found in detail in \cite[Ch. 10]{AH} and \cite{cia}.

 For a given interval matrix $\mc=[\mcs_{ij}]$
let us define corresponding point matrices, the {\it midpoint
matrix} $m(\mc):=[m(\mcs_{ij})],$  the {\it width matrix} $d(\mc):=[d(\mcs_{ij})],$ and the {\it
absolute value matrix} $|\mc|:=[|\mcs_{ij}|],$
 as follows:
  $$m(\mcs_{ij})=\tfrac12(\underline{c}_{\,ij}+\overline{c}_{ij}),
\quad d(\mcs_{ij}):=\overline{c}_{\,ij}-\underline{c}_{\,ij}, \quad
|\mcs_{ij}|=\max\{|\overline{c}_{ij}|,|\underline{c}_{\,ij}|\}.
$$ If $\mc=\dc=[c_{ij}] $ is a point matrix, then it is obvious
$$m(\dc)=[c_{ij}],\quad d(\dc)=[0]\ \mbox{\rm (null-matrix)},\quad
|\dc|=[|c_{ij}|].
$$

We start with the following
 matrix identity for an
 $n\times n$ matrix $\dq$ and the unity matrix $\di,$
 $$
 (\di-\dq)(\di+\dq+\cdots+\dq^{r-2})=\di-\dq^{r-1}.
 $$
  Hence, setting $\dq=\da\dh,$ where $\dh$ is an $n\times n$ matrix, the following identity is obtained:
 \bb \dh\sum_{\la=0}^{r-2}(\di-\da\dh)^{\la}=\da^{-1}-\da^{-1}
 (\di-\da\dh)^{r-1}.\label{i1}
 \ee
 From (\ref{i1}) there follows
 \bb \da^{-1}=\dh\sum_{\la=0}^{r-2}(\di-\da\dh)^{\la}+\da^{-1}(\di-\da\dh)^{r-1}.\label{i2}
 \ee
This relation will be used for the construction of  interval matrix iterations.

Let $\mx_0$ be an  $n\times n$ interval matrix such that $\da^{-1}\in
\mx_0,$ and let the matrix $\dh$ in (\ref{i2}) be defined by
$\dh=m(\mx_0).$ Then we obtain from (\ref{i2}) using inclusion property
 \bb \da^{-1}\in \mx_1:=m(\mx_0)\sum_{\la=0}^{r-2}
 \Bigl(\di-\da m(\mx_0)\Bigr)^{\la}+
 \mx_0\Bigl(\di-\da m(\mx_0)\Bigr)^{r-1}.\label{i3}
 \ee

For simplicity, let us introduce $\dr_k=\di-\da m(\mx_k).$ Combining  (\ref{i2}) and (\ref{i3}), it is easily to prove by the
set property and mathematical induction that the following is
valid for an arbitrary $k\ge 0:$ \bb \da^{-1}\in
\mx_{k+1}:=m(\mx_k)\sum_{\la=0}^{r-2}
\dr_k^{\la}+
 \mx_k\dr_k^{r-1}.\label{i4}
 \ee
 In regard to this property, the following iterative
 process for finding an inclusion  matrix for $\da^{-1}$ can
 be stated in a Horner scheme fashion
 \bb \left\{\begin{array}{l}
 \my_{k}=m(\mx_k)\Bigl(\di\underbrace{+\dr_k(\di+\dr_k(\di+\cdots+
 \dr_k(\di+}_{r-2\ \mbox{\rm times}}\dr_k)\cdots \Bigr)

 +
 \mx_k \dr_k^{r-1},\\[10pt]
 \mx_{k+1}=\my_k\cap \mx_k,
 \end{array}\right.\quad (k=0,1,\ldots).\label{i5}
 \ee

 The iterative method (\ref{i5}) was considered in
 detail in the book \cite{AH} by Alefeld and Herzberger.
 As shown in \cite[Ch. 18]{AH}, the most efficient method from the class (\ref{i5}) of hyper-power methods is obtained for $r=3$ and reads
    \bb \left\{\begin{array}{l}
 \my_{k}=m(\mx^{(k)})+m(\mx^{(k)})\dr_k+
  \mx_k \dr_k^2,\\[10pt]
 \mx_{k+1}=\my_k\cap \mx_k,
 \end{array}\right.\quad (k=0,1,\ldots).\label{i6a}
 \ee

The properties of the iterative interval method (\ref{i5}) are
 given in the following theorem proved in \cite[Theorem 2, Ch.
 18]{AH}, where $\rho(M)$ denotes the spectral radius of a matrix $M.$


\begin{thm} \fontsize{11}{12}\selectfont Let $\da$ be a nonsingular $n\times n$ matrix
 and $\mx_0$ an $n\times n$ interval matrix such that
 $\da^{-1}\in \mx_0.$ Then

\begin{itemize}
\itemsep0pt
\item[(a)]
 each inclusion matrix $\mx_k,$ calculated by
 $(\ref{i5})$, contains $A^{-1};$
\item[(b)] If $\rho(|\di-\da\dx|)<1$ for every $\dx \in \mx_0,$ then
the sequence $\{\mx_k\}_{k\ge 0}$ converges to
$A^{-1};$
 \item[(c)] using a matrix norm $\|\cdot\|$ the sequence
 $\{d(\mx_k)\}_{k\ge 0}$ satisfies
  $$
  \|d(\mx_{k+1})\|\le \gamma \|d(\mx_k)\|^r,\quad \gamma\ge 0,
  $$
  that is, the $R$-order of convergence of the method $(\ref{i5})$ is at least $r.$
  \end{itemize} 
\label{thm:i1}
\end{thm}


 Using the iterative formula (\ref{i5}) in the Horner form for $r=6,$ we obtain the following iterative method for the inclusion of the inverse matrix:
    \bb \aligned
    \dr_k&=\di -\da\odot m(\mx_k),\\
    \ds_k&=\dr_k\odot \dr_k,\\
 \dm_k&=\di+\dr_k\odot(\di+\dr_k\odot(\di+\dr_k\odot
 (\di+\dr_k))),\\
 \my_k&= m(\mx_k)\odot \dm_k+\mx_k \otimes (\ds_k\odot \ds_k\odot \dr_k),\\
 \mx_{k+1}&=\my_k\cap \mx_k,
 \endaligned
 \qquad (k=0,1,\ldots).
 \label{i7}
 \ee
 The method \eqref{i7} is a particular case of the general matrix iteration \eqref{i5}.
 According to Theorem \ref{thm:i1}, the method \eqref{i7} has order six
 and requires 8 multiplication of point
 matrices (denoted by $\odot$) and one multiplication
 of interval matrix by  point matrix (denoted by $\otimes$).

\section{\bsfs New inclusion method of high efficiency}

 In what follows we are going to show that the computational cost of the interval method (\ref{i7})
 can be reduced using the identity
  \bb
  x^4+x^3+x^2+x+1=x^2(x^2+x+1)+x+1
  \label{rel}
  \ee
  and the corresponding matrix relation. Having in mind (\ref{rel}) we rewrite (\ref{i7}) and
 state the following algorithm in interval arithmetic for bounding the inverse matrix:
 \bb
\aligned
\dr_k&=\di-\da\odot m(\mx_k),\\
\ds_k&=\dr_k\odot \dr_k,\\
\dt_k&=\ds_k\odot \ds_k\odot \dr_k,\\
\dm_k&=\di+\dr_k+\ds_k\odot (\di+\dr_k+\ds_k),\\
\my_k&=m(\mx_k)\odot \dm_k+\mx_k\otimes \dt_k,\\
\mx_{k+1}&=\my_k\cap \mx_k,
\endaligned
 \qquad (k=0,1,\ldots).
\label{i8} \ee Compared with the method (\ref{i7}), the iterative
scheme (\ref{i8}) requires  6 multiplications of point matrices
(thus, two matrix multiplications less) and still preserves the
order six. The above consideration can be summarized in the following
theorem.

\begin{thm} \fontsize{11}{12}\selectfont Let $\da$ be a nonsingular $n\times n$ matrix
 and $\mx_0$ an $n\times n$ interval matrix such that
 $\da^{-1}\in \mx_0.$ Then

 \smallskip

\begin{itemize} 
\itemsep0pt
\item[(a)]
 each inclusion matrix $\mx_k,$ calculated by
 $(\ref{i8})$, contains $\da^{-1};$
\item[(b)]  if $\rho(|\di-\da \dx|)<1$ holds for all $\dx\in \mx_0,$ then the sequence $\{\mx_k\}_{k\ge 0}$
   converges toward $\da^{-1}$;
 \item[(c)] using a matrix norm $\|\cdot\|$ the sequence
 $\{d(\mx_k)\}_{k\ge 0}$ satisfies
  $$
  \|d(\mx_{k+1})\|\le \gamma \|d(\mx_k)\|^6,\quad \gamma\ge 0,
  $$
  that is, the $R$-order of convergence of the
  method $(\ref{i8})$ is at least $6.$
  \end{itemize}
\label{thm:i2}
\end{thm}

  \noindent  Theorem \ref{thm:i2} can be proved in a similar way
  as Theorems 1
  and 2 in \cite[Ch. 18]{AH} so that we omit the proof.


\begin{rem} \it \fontsize{11}{12}\selectfont
Zhang, Cai and Wei have proved in \cite[Theorem 3.3]{Interval} that, under the additional condition $(m(\mx_0)=A^TBA^T$ for some matrix $B\in \R^{m\times m})$, the iterative method \eqref{i5} (and specially \eqref{i7}) is also convergent in the case of full-rank $m\times n$ matrix $\da$. In such a case, it converges to the Moore-Penrose inverse $\da^\dagger$ of $\da$. In a similar way, the same can be proved for the
method \eqref{i8}.

\end{rem}

Executing iterative interval processes in general, one of the most important
but also very difficult task is to find a good initial interval
(real interval, complex interval, interval matrix, etc.) that
contains the sought result. Similar situation appears in bounding
the inverse matrix. We present here an efficient method for
construction an initial matrix $\mx_0$ that contains the inverse
matrix $\da^{-1}.$

 Let $\dx\in \mx_0$
and let us assume that the matrix $\da$ can be represented as
 \bb  \da=\di-\dy,\quad \mbox{where a newly introduced matrix $\dy$ satisfies}\quad
\|\dy\|<1.\label{i9}
  \ee
   It has been shown in  \cite[Ch. 18]{AH} that the inequality
 $$
 \|\dx\|\le a:=\frac{1}{1-\|\dy\|}
 $$
 holds.
 If we use either the row-sum or the column-sum norm, then we find that
  $
  -a\le x_{ij}\le a,\ (1\le i,j\le n)
  $
  holds for all the elements of $\dx=[x_{ij}].$ For the matrix
  $\mx_0=\bigl[ X_{ij}^{(0)}\bigr]$ with
  interval coefficients
   \bb
   X_{ij}^{(0)}=\left\{\begin{array}{ll}
   [-a,a] &  \mbox{\rm for}\ i\ne j\\[2pt]
   [-a,2+a] & \mbox{\rm for}\ i=j,\end{array}\right.\label{i10}
   \ee
   we have $\da^{-1}\in \mx_0$ and $m(\mx_0)=\di$ (see \cite{AH}).
   If the condition (\ref{i9}) is not satisfied,
   then it is effectively to normalize
   the matrix $\da$ before running the
   iterative process, say, to deal with the matrices $\da/\|\da\|$ or $\da/\|\da\|^2.$

   Having in mind the described procedure of choosing
   initial inclusion matrix $\mx_0,$
   applying point matrix iterations
   it is convenient to take $\dx_0=m(\mx_0)=\di.$ Such choice have already
   applied in stating the iterative interval methods (\ref{i7}) and (\ref{i8}).

   \section{\bsfs Computational aspects}

   Let us compare computational efficiency of the hybrid methods (\ref{i7}) and (\ref{i8}). As proved in \cite[Ch. 6]{springer}, CPU (central processor unit) time necessary for executing an iterative method $(IM)$ can be suitably expressed in a pretty manner in the form
    \bb
      CPU_{\footnotesize (IM)}=h\log q \cdot \frac{\theta(IM)}
      {\log r(IM)}.\label{e1}
      \ee
      Here $r(IM)$ is the convergence order, $\theta(IM)$ is computational cost of the iterative method $(IM)$ per iteration, $q$ is the number od significant decimal digits (for example $q=15$ or $16$ for double precision arithmetic) and $h$ is a constant that depends on hardware characteristics of the employed digital computer. Assuming that the considered methods are implemented on the same computer, according to (\ref{e1}) the comparison of two methods $(M_1)$ and $(M_2)$ is carried out  by the {\it efficiency ratio}
       \bb
       ER_{\small M_1/ M_2}(n)=\frac{CPU_{\footnotesize (M_1)}}
       {CPU_{\footnotesize (M_2)}}=\frac{\log r(M_2)}{\log r(M_2)}
      \cdot \frac{\theta(M_1)}{\theta(M_2)}.\label{e2}
       \ee

Calculating the computational cost $\theta,$ it is necessary to deal with  the number of arithmetic operations per
iteration taken with certain {\it weights} depending on the
 execution times of operations.
We  assume
that floating-point number representation is used, with a binary
fraction of $b$ bits, meaning that we deal with ``{\it
precision} $b$" numbers, giving results with  a relative error of
approximately $2^{-b}.$ Following results given in \cite{brent},
 the execution time $t_b(A)$
of addition (subtraction) is $\od(b),$ where $\od$ is the Landau symbol. Using
Sch\"onhage-Strassen multiplication (see \cite{brent}),
 often implemented in multi-precision libraries
(in the computer algebra systems  {\it Mathematica, Maple, Magma},
for instance),
we have  $t_b(M)=\od\bigl(b\log b\, \log (\log b)\bigr).$ For comparison purpose, we chose the weights $w_{a}$ and $w_m$
proportional to $t_b(A)$ and $t_b(M),$ respectively for double precision arithmetic ($b=64$ bits) and quadruple-precision  arithmetic ($b=128$ bits).

In particular cases, assuming that multiplication  of two scalar $n\times n$ matrices requires $n^2(n-1)$ additions and $n^3$ multiplications and adding combined costs in the iterative formulae (\ref{i7}) and (\ref{i8}), for the hybrid method (\ref{i7}) and (\ref{i8}) we have $r(\ref{i7})=r(\ref{i8})=6$ and, approximately,
$$
\theta(\ref{i7})=(9n^3-3n^2)b+10n^3b\log b \log(\log b),\quad
\theta(\ref{i8})=(7n^3-n^2)b+8n^3b\log b \log(\log b).
$$
In view of this, by (\ref{e2}) we determine the {efficiency ratio}
 $$
 ER_{(\ref{i7})/(\ref{i8})}(n)=\frac{9-3/n+10\log b \log(\log b)}
 {7-1/n+8\log b \log(\log b)}.
$$

The graph of the function $ER_{(\ref{i7})/(\ref{i8})}(n)$ for $n\in [2,40]$ is shown in Figure 1. From this graph we note that the values of
$ER(n)$ are grouped about the value 1.25 for $n$ in a wide range. This means that the new method (\ref{i8}) consumes about 25\% less CPU time than
the Horner-fashion method (\ref{i7}).

\cl{\includegraphics[height=6cm]{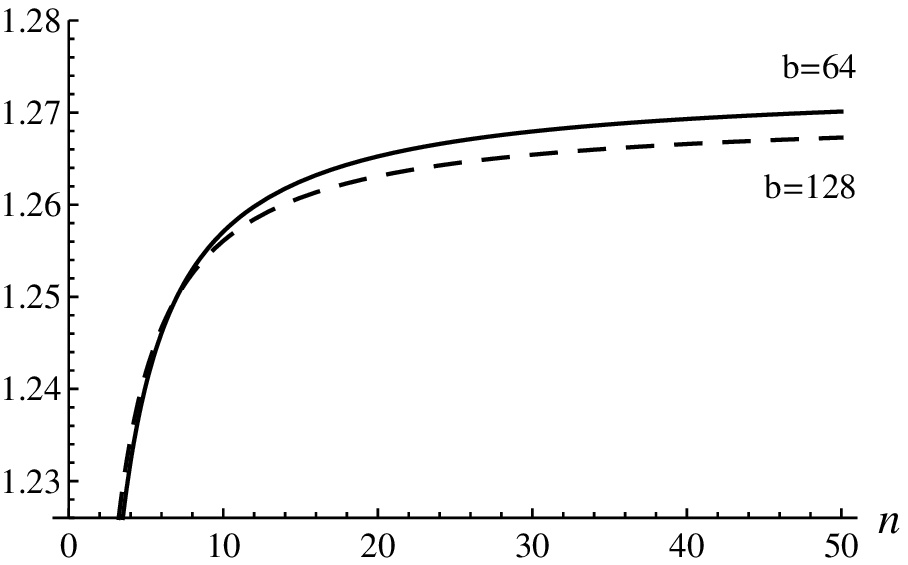}}



\centerline{\small Figure 1: The ratio of CPU times for two different precisions of arithmetical processors}

\bigskip

 A very similar graph is obtained for a lot of computing machines. For example, for double precision arithmetic  and quadruple precision arithmetic (corresponding approximately to $b=64$ and $b=128,$ respectively) for the processor Pentium M 2.8 GHz (Fedora core 3) the values of $ER(n)$ are very close to 1.25 almost independently on the dimension of matrix $n.$ In addition, we find $ER_{(\ref{i5})/(\ref{i8})}(n)>1$ for every $r\ne 3$ and close to 1 for $r=3.$


   The convergence behavior of the iterative interval method (\ref{i8}), together with the choice of initial inclusion matrix $\mx_0,$ will be demonstrated by one  simple example. We emphasize that  the interval method (\ref{i7}) produces the same inclusion matrix, which is obvious since the corresponding iterative formulae are, actually, identical but arranged in different forms. However,
   as mentioned above, the inclusion method (\ref{i8}) has
 lower computational cost than (\ref{i7}).

   \bigskip

   \noindent {\bf Example 1.} We wish to find the inclusion matrix for the inverse of the matrix
     $$
     \da=\left[\begin{array} {rr}
     \frac{9}{10} & \frac{1}{5}\\[6pt] -\frac{3}{10} & \frac{4}{5}
     \end{array}\right].
     $$ Note the the inverse matrix $\da^{-1}$ is
          $$
     \da^{-1}=\left[\begin{array} {rr}
     \frac{40}{39}  & -\frac{10}{39}\\[6pt] \frac{5}{13} & \frac{15}{13}
     \end{array}\right]=\left[\begin{array} {rr}
     1.0\overline{256410} & -0.\overline{256410} \\[2pt]
     0.\overline{384615} & 1.\overline{153846}
     \end{array}\right].
     $$
 The overlined set of digits indicates that this set of digits repeats periodically.

     First we determine
      $$
      \dy=\di-\da=\left[\begin{array} {rr}
     0.1 & -0.2\\ 0.3 & 0.2
     \end{array}\right] \ \ \mbox{\rm with} \ \ \|\dy\|_2=0.424264\ \ \mbox{\rm and}\ \ a=\frac{1}{1-\|\dy\|_2}=1.73691.
     $$
   According to (\ref{i10}) we form the initial inclusion matrix
        $$
     \mx_0=\left[\begin{array} {ll}
     [-1.73691,3.73691] & [-1.73691,1.73691]\\[2pt]
  [-1.73691,1.73691] & [-1.73691, 3.73691]
     \end{array}\right].
     $$
     Note that the widths of intervals which present the coefficients
     of the initial inclusion matrix $\mx_0$ are rather large.
We have applied two iterations of (\ref{i8}) and obtained the following midpoint matrices (approximations to $\da^{-1}$) and the width matrices that give the upper error bounds of $\mx_k.$

\smallskip
\hskip1.5cm$\underline{k=1}$
  \ba
   &&m(\mx_1)=\left[\begin{array} {rr}
    1.025\ldots  & -0.256\ldots\\[2pt]
 0.384\ldots  & 1.153\ldots
     \end{array}\right],\quad d(\mx_1)=\left[\begin{array} {ll}
  1.27\times 10^{-2}  & 8.68\times 10^{-3}\\[2pt]
 1.51\times 10^{-2} & 6.356\times 10^{-3}
 \end{array}\right].
             \ea

\hskip1.5cm$\underline{k=2}$
\ba
      && m(\mx_2)=\left[\begin{array} {rr}
    1.0256410256410256\ldots  & -0.256410256410256\ldots\\[2pt]
 0.3846153846153846\ldots  &  1.153846153846153\ldots
     \end{array}\right], \\ && 
     d(\mx_2)=\left[\begin{array}{ll}
  6.33\times 10^{-19}  & 4.19\times 10^{-19}\\[2pt]
 5.99\times 10^{-19} & 4.54\times 10^{-19}
     \end{array}\right].
        \ea
        All displayed decimal digits of $m(\mx_1)$ and $m(\mx_2)$ are correct.
       The third iteration produces the
       width matrix  $d(\mx_3)$ with elements
       in the form of real intervals with widths of order $10^{-99}.$
       We have not listed $m(\mx_3)$ and $d(\mx_3)$ to save the space.

We have also tested the interval method (\ref{i6a}) possessing the highest efficiency among hyper-power methods. Starting with the same initial matrix $\mx_0$ as above, we obtained the following outcomes:

\smallskip
\hskip1.5cm$\underline{k=1}$
  \ba
   &&m(\mx_1)=\left[\begin{array} {rr}
    1.05  & -0.26\\[2pt]
 0.39  & 1.18
     \end{array}\right],\quad d(\mx_1)=\left[\begin{array} {ll}
  0.586  & 0.398\\[2pt]
 0.666 & 0.318
 \end{array}\right].
             \ea

\hskip1.5cm$\underline{k=2}$
\ba
      && m(\mx_2)=\left[\begin{array} {rr}
    1.0256\ldots  & -0.2564\ldots\\[2pt]
 0.3846\ldots  &  1.1538\ldots
     \end{array}\right],\quad d(\mx_2)=\left[\begin{array}{ll}
  3.60\times 10^{-4}  & 2.43\times 10^{-4}\\[2pt]
 3.91\times 10^{-4} & 2.12\times 10^{-4}
     \end{array}\right].
        \ea

 The method (\ref{i8}) produced considerably higher accuracy than (\ref{i6a}) using only two iterations so that its application is justified in this case.
 Furthermore, since $
  ER_{(\ref{i6a})/(\ref{i8})}(n)
  $
 is  close to 1,
 which of these  two methods will be chosen depends of the nature of solved problem, specific requirements and available hardware and software (precision of employed computer). For instance, the proposed method (\ref{i8}) is more convenient when a high accuracy is requested in a few iterations, as in the presented example.

\end{document}